\documentclass[12pt]{amsart}%
\usepackage{amsfonts}
\usepackage{txfonts}
\usepackage{mathtools}
\usepackage{amsmath}
\usepackage{amssymb}
\usepackage{amsthm}
\usepackage{newlfont}
\usepackage{graphicx}%
\providecommand{\U}[1]{\protect\rule{.1in}{.1in}}
\newtheorem{theorem}{Theorem}[section]
\theoremstyle{plain}

\newtheorem{corollary}{Corollary}[section]

\newtheorem{remark}{Remark}[section]

\numberwithin{equation}{section}

\begin{document}

\title[Hamiltonian stationary cone]{Hamiltonian stationary cones with isotropic links }

\dedicatory{In memory of Professor Wei-Yue Ding}

\author{Jingyi CHEN}
\address{Department of Mathematics\\
University of British Columbia\\
Vancouver, B.C., V6T 1Z2\\
Canada}
\email{jychen@math.ubc.ca}
\author{Yu YUAN}
\address{Department of Mathematics, Box 354350\\
University of Washington\\
Seattle, WA 98195\\
USA}
\email{yuan@math.washington.edu}
\thanks{The first authors is partially supported by NSERC, and the second author is
partially supported by an NSF grant.}
\date{\today}

\begin{abstract}
We show that any closed oriented immersed Hamiltonian stationary isotropic
surface $\Sigma$ with genus $g_{\Sigma}$ in $S^{5}\subset\mathbb{C}^{3}$ is
(1) Legendrian and minimal if $g_{\Sigma}=0$; (2) either Legendrian or with
exactly $2g_{\Sigma}-2$ Legendrian points if $g_{\Sigma}\geq1.$ In general,
every compact oriented immersed isotropic submanifold $L^{n-1}\subset
S^{2n-1}\subset\mathbb{C}^{n}$ such that the cone $C\left(  L^{n-1}\right)  $
is Hamiltonian stationary must be Legendrian and minimal if its first Betti
number is zero. Corresponding results for non-orientable links are also provided.
\end{abstract}

\maketitle

\section{\bigskip Introduction}

In this note we study the problem of when a Hamiltonian stationary cone $C(L)$
with isotropic link $L$ on $S^{2n-1}\ $ in $\mathbb{C}^{n}$ becomes special
Lagrangian. A submanifold $M\subset\mathbb{C}^{n}$, not necessarily a Lagrangian submanifold, is {\it Hamiltonian stationary}
means%
\[
\operatorname{div}_{M}\left(  JH\right)  =0,
\]
where $J$ is the complex structure in $\mathbb{C}^{n}$ and $H$ is the mean curvature vector of $M$ in ${\mathbb C}^n$. In fact this is the
variational equation of volume of $M$, when one makes an arbitrary deformation
$J\nabla_{M}\varphi$ with $\varphi\in C_{0}^{\infty}\left(  N\right)  $ for
$N$:
\[
\int_{M}\left\langle H,J\nabla_{M}\varphi\right\rangle =\int_{M}%
\varphi\operatorname{div}_{M}\left(  JH\right)  -\operatorname{div}_{M}\left(
\varphi JH\right)  =\int_{M}\varphi\operatorname{div}_{M}\left(  JH\right)  .
\]
The notion of Hamiltonian stationary Lagrangian submanifolds in a K\"ahler manifold was introduced by Oh \cite{Oh} as critical points of the volume functional 
under Hamiltonian variations (known to A. Weinstein as noted in \cite{Oh}).  Chen-Morvan \cite{CM} generalized it to the isotropic deformations. 

As in Harvey-Lawson \cite{HL}, a submanifold $M$
in $\mathbb{C}^{n}$ is {\it isotropic} at $p\in M$ if 
\[
J\left(  T_pM\right)  \perp T_pM,
\]
and it is isotropic if it is isotropic for every $p$. A submanifold $M$ being isotropic is equivalent to that 
 the standard symplectic 2-form on $\mathbb{R}^{2n}$ vanishes on $M$. The
dimension of an isotropic submanifold is at most $n$, the half real dimension of ${\mathbb C}^n$, and
when it is $n$ the submanifold is Lagrangian. 

For an immersed $(n-1)$-dimensional submanifold $L$ in the unit sphere $S^{2n-1}$, let $u:L\rightarrow S^{2n-1}$ be
the restriction of the coordinate functions in ${\mathbb{R}}^{2n}$ to $L$. A
point $u\in L$ is {\it Legendrian} if $T_{u}L$ is isotropic in $\mathbb{R}^{2n}$
and
\[
J\left(  T_{u}L\right)  \perp u.
\]
$L$ is Legendrian if all the points $u$ are Legendrian. This is equivalent to
that $L$ is an $\left(  n-1\right)  $-dimensional integral submanifold of the
standard contact distribution on $S^{2n-1}$. The cone $C(L)= \{r x:r\geq 0, x\in L\}$ is said to have {\it link} $L$. 
In this article, all links $L$ are assumed to be connected, and we shall use $\Sigma$ for 2-dimensional link $L$.

The Hamiltonian stationary condition is a third order constrain on the submanifold $M$, as seen when $M$ is locally written as a graph over 
its tangent space at a point. The minimal submanifolds, a second order constrain on the local graphical representation of $M$, are 
automatically Hamiltonian stationary. We are particularly interested in the case that $M$ is a Lagrangian submanifold. The existence of  (many) compact Hamiltonian stationary Lagrangian submanifolds in ${\mathbb C}^n$ versus the non-existence of compact minimal submanifolds makes the study of Hamiltonian stationary ones interesting. In this note, we shall not be concerned with the existence of Hamiltonian stationary ones; instead, we shall concentrate on the rigidity property, namely, when the Hamiltonian stationary ones reduce to special Lagrangians, in the case when the submanifold is a cone over a spherical link in ${\mathbb C}^n$. 

A well-known fact about a link $L^m\subset S^n$ and the cone $C(L)$ over it is that $L$ is minimal in $S^n$ if and only if $C(L)\setminus\{0\}$ is minimal in ${\mathbb R}^{n+1}$. When $C(L)$ is Hamiltonian stationary and isotropic, possibly away from the cone vertex $0\in{\mathbb R}^{2n}$, we observe that the Hamiltonian stationary equation for $C(L)$ splits into two equations:
$$
 \operatorname{div}_{L}\left(JH_L\right)  =0       
$$
i.e. the link $L$ is Hamiltonian stationary in ${\mathbb R}^{2n}$ as well, and 
$$
\langle JH_L, u \rangle =0,
$$
where $H_L$ is the mean curvature vector of $L$ in ${\mathbb R}^{2n}$ and $u$ is the position vector of $L$. Moreover, if the link $L$ is isotropic in ${\mathbb C}^n$, then 
$$
  \operatorname{div}_{L}\left(J\bar{H}_{L}\right)=0
$$  
where $\bar{H}_L =H_L-mu$ is the mean curvature vector of $L$ in $\mathbb S^{2n-1}$; in fact,
$$
 \operatorname{div}_{L}(Ju) =\sum_{i=1}^m\langle D_{E_i} (Ju), E_i\rangle
 =\sum^m_{i=1}\langle JD_{E_i}u,E_i\rangle = 0
 $$
as $D_{E_i}u$ is tangent to $L$, where $D$ is the derivative in ${\mathbb R}^{2n}$ and $\{E_1, \cdots, E_m\}$ is an orthonormal local frame on $TL$.

Our observation is that the rigidity statements in \cite{CY}  for minimal links  generalize to the Hamiltonian stationary setting.

\begin{theorem}
\label{thm-Hopf}Let $\Sigma$ be a closed oriented immersed isotropic surface
with genus $g_{\Sigma}$ in $S^{5}\subset\mathbb{C}^{3}$ such that the cone
$C\left(\Sigma\right)$ is Hamiltonian stationary away from its vertex. Then 

(1) if $g_{\Sigma}=0$,
$\Sigma$ is Legendrian and minimal (in fact, totally geodesic); 

(2) if $g_{\Sigma}\geq1$, $\Sigma$ is
either Legendrian or has exactly $2g_{\Sigma}-2$ Legendrian points counting
the multiplicity.
\end{theorem}

It is known that the immersed minimal Legendrian sphere ($g_{\Sigma}=0$) must
be a great two sphere in $S^{5}$ (cf. [H, Theorem 2.7]). Simple isotropic tori
($g_{\Sigma}=1$) can be constructed so that the Hamiltonian stationary cone
$C\left(  \Sigma\right)  $ is nowhere Lagrangian. A family of Hamiltonian
stationary (non-minimal) Lagrangian cones $C\left(  \Sigma\right)  $ with
$g_{\Sigma}=1$ are presented by Iriyeh in [I]. Bryant's classification [Br,
p.269] of minimal surfaces with constant curvature in spheres provides
examples of flat Legendrian minimal tori as well as flat non-Legendrian
isotropic minimal tori ($g_{\Sigma}=1$) in $S^{5}$. The constructions of
Haskins [H] and Haskins-Kapouleas [HK] show that there are infinitely many
immersed (embedded if $g_{\Sigma}=1$) minimal Legendrian surfaces for each odd
genus in $S^{5}$.

In general dimensions and co-dimensions, we have

\begin{theorem}
\label{thm-Hodge}Let $L^{m}$ be a compact isotropic immersed oriented
submanifold in the unit sphere $S^{2n-1}\subset\mathbb{C}^{n}$ such that the
cone $C\left(  L^m\right)  $ is Hamiltonian stationary away from its vertex. Suppose that the first Betti
number of $L^{m}\ $is $0$.  Then, away from its vertex,

 (1) when $m$ is the top dimension $n-1$, the cone $C(L^{n-1})$ is Lagrangian and minimal (or equivalently $L^{n-1}$
is Legendrian and minimal); 

(2) for $m<n-1$, the cone $C(L^{m})$ is isotropic, and if the differential 1-form $\langle JH_{C(L^m)},\cdot\rangle$ is closed then the mean curvature $H_{C(L^m)}$ of $C(L^m)$ vanishes on the normal subbundle $JTC(L^m)$.
\end{theorem}



We make two remarks when the dimension $m$ of the link is two. First, Theorem \ref{thm-Hodge} also
implies Theorem \ref{thm-Hopf} (1). Second, if the first Betti number of $L^2$ is not zero
($g_{L^{2}}>0$) and $L$ is isotropically immersed in $S^{2n-1}$, with $2n-1\geq 5$, and $C(L)$ is Hamiltonian stationary away from its cone vertex, the same argument as in the proof of Theorem 1.1 leads to the same conclusion as in part (2) of Theorem 1.1 that the
cone $C(L^{2})$ is isotropic either everywhere or along exactly $2g_{L^{2}
}-2=-\chi (  L^{2} )$ lines. 

Theorem \ref{thm-Hodge} and Theorem \ref{thm-Hopf} (except
the totally geodesic part) remain 
valid for non-orientable links (note that $\chi\left(  \Sigma\right)
=2-g_{\Sigma}$ for a compact non-orientable surface $\Sigma$), see Remark 2.2
and Remark 3.1. The non-orientable version of Theorem
\ref{thm-Hodge} implies that one cannot immerse a compact non-orientable
$L^{n-1}$ with zero first Betti number Hamiltonian stationarily and
isotropically into $S^{2n-1}\subset\mathbb{C}^{n}.$ Otherwise, the cone
$C(L^{n-1})$ would be a special Lagrangian cone, then $C(L^{n-1})$ would be
orientable, and $L^{n-1}$ would also be orientable. In particular, there exists no isotropic Hamiltonian stationary
immersion of a real projective sphere $\mathbb{R}P^{2}$ into $S^{5}
\subset\mathbb{C}^{3}.$

It is interesting to find out whether there exists an isotropic Hamiltonian stationary
surface in $S^{5}$ with exactly $2g_{\Sigma}-2$ Legendrian points for
$g_{\Sigma}>1$.


\section{ Hopf differentials and Proof of theorem \ref{thm-Hopf}}

To measure how far the isotropic $\Sigma$ is away from being Legendrian, or
the deviation of the corresponding cone from being Lagrangian, we project $Ju$
onto the tangent space of $\Sigma$ in $\mathbb{C}^{3},$ where $J$ is the
complex structure in $\mathbb{C}^{3}.$ Denote the length of the projection by
\[
f=\left\vert \Pr Ju\right\vert ^{2}.
\]
To compute the length, we need some preparation. Locally, take an isothermal
coordinate system $\left(  t^{1},t^{2}\right)  $ on the isotropic surface
\[
u:\Sigma\rightarrow S^{5}\subset\mathbb{C}^{3}.
\]
Set the complex variable
\[
z=t^{1}+\sqrt{-1}t^{2}.
\]
Then the induced metric has the local expression with the conformal factor
$\varphi$
\[
g=\varphi^{2}\left[  \left(  dt^{1}\right)  ^{2}+\left(  dt^{2}\right)
^{2}\right]  =\varphi^{2}dzd\bar{z}.
\]
We project $Ju$ to each of the orthonormal basis $\varphi^{-1}u_{1}%
,\varphi^{-1}u_{2}$ with $u_{i}=\partial u/\partial t^{i}$. Then the sum of
the squares of each projection is%

\[
f=\frac{\left\vert \left\langle Ju,u_{1}\right\rangle \right\vert
^{2}+\left\vert \left\langle Ju,u_{2}\right\rangle \right\vert ^{2}}%
{\varphi^{2}}=\frac{4\left\vert \left\langle Ju,u_{z}\right\rangle \right\vert
^{2}}{\varphi^{2}}%
\]
where $u_{z}=\partial u/\partial z$ and $\langle\cdot,\cdot\rangle$ is the
Euclidean inner product on $\mathbb{R}^{6}$, and in particular $0\leq f\leq1$.
In fact, $f$ is the square of the norm of the symplectic form $\omega$ in
$\mathbb{C}^{3}$ restricted on the cone $C(\Sigma)$ with link $\Sigma$:
\[
\omega|_{C(\Sigma)}\wedge\ast\,\omega|_{C(\Sigma)}=f\,\cdot
\,\hbox{volume form of $C(\Sigma)$}.
\]

The Hamiltonian stationary condition for the cone 
$
C\left(  \Sigma\right)=ru\left(  t^{1},t^{2}\right)  $ is
\begin{eqnarray*}
0  &  =&\operatorname{div}_{C\left(  \Sigma\right)  }\left(  JH_{C\left(
\Sigma\right)  }\right) \\
 &=&\left\langle \partial_{r}\left(  JH_{C\left(
\Sigma\right)  }\right)  ,\partial_{r}\right\rangle +\frac{1}{r^{2}%
}\operatorname{div}_{\Sigma}\left(  J\frac{1}{r}H_{\Sigma}\right) \\
&  =&-\frac{1}{r^{2}}\left\langle JH,u\right\rangle +\frac{1}{r^{3}%
}\operatorname{div}_{\Sigma}\left(  JH_{\Sigma}\right)  .
\end{eqnarray*}
It follows that%
\[
\operatorname{div}_{\Sigma}\left(  JH_{\Sigma}\right)  =0
\]
and%
\[
0=\left\langle JH,u\right\rangle =-\left\langle \frac{4}{\varphi^{2}}%
u_{z\bar{z}},Ju\right\rangle .
\]
Coupled with the isotropy condition
\[
\,\left\langle Ju_{i},u_{j}\right\rangle =0
\]
we have the holomorphic condition
\[
\left\langle Ju,u_{z}\right\rangle _{\bar{z}}=\left\langle Ju_{_{\bar{z}}%
},u_{z}\right\rangle +\left\langle Ju,u_{z\bar{z}}\right\rangle =\left\langle
Ju,-\frac{\varphi^{2}}{2}u\right\rangle =0.
\]

The induced metric $g$ yields a compatible conformal structure on the oriented
surface $\Sigma$ which makes $\Sigma$ a Riemann surface. We shall consider two
cases according to the genus $g_{\Sigma}$.

(1) $g_{\Sigma}=0$: By the uniformization theorem for Riemann surfaces (cf.
[AS, p.125, p.181]), there exists a holomorphic covering map
\[
\Phi:\left(  S^{2},g_{\text{canonical}}\right)  \rightarrow\left(
\Sigma,g\right)
\]
or locally
\[
\Phi:\left(  \mathbb{C}^{1},\frac{1}{\left(  1+\left\vert w\right\vert
^{2}\right)  ^{2}}dwd\bar{w}\right)  \rightarrow\left(  \Sigma,g\right)  .
\]
For $z=\Phi(w)$ one has
\[
\frac{1}{\left(  1+\left\vert w\right\vert ^{2}\right)  ^{2}}dwd\bar{w}%
=\Phi^{\ast}\left(  \mathit{\psi}^{2}g\right)  =\Phi^{\ast}\left(
\mathit{\psi}^{2}\varphi^{2}dzd\bar{z}\right)  =\mathit{\psi}^{2}\varphi
^{2}\left\vert z_{w}\right\vert ^{2}dwd\bar{w},
\]
where $\psi$ is a positive (real analytic) function on $\Sigma.$ In
particular
\[
\left\vert z_{w}\right\vert ^{2}=\frac{1}{\psi^{2}\varphi^{2}\left(
1+\left\vert w\right\vert ^{2}\right)  ^{2}}.
\]
Note that
\[
\left\langle Ju,u_{w}\right\rangle =\left\langle Ju,u_{z}\right\rangle
z_{w}=\left\langle Ju,u_{z}\right\rangle \frac{1}{w_{z}}%
\]
is a holomorphic function of $z$, in turn it is a holomorphic function of $w$.
Also $\left\langle Ju,u_{w}\right\rangle $ is bounded, approaching $0$ as $w$
goes to $\infty,$ because%
\[
\left\vert \left\langle Ju,u_{w}\right\rangle \right\vert ^{2}=\frac
{\left\vert \left\langle Ju,u_{z}\right\rangle \right\vert ^{2}}{\varphi^{2}%
}\frac{1}{\psi^{2}\left(  1+\left\vert w\right\vert ^{2}\right)  ^{2}}.
\]
So $\left\langle Ju,u_{w}\right\rangle \equiv0.$ Therefore $f\equiv0$ and
$\Sigma$ is Legendrian. We conclude that $C(\Sigma)\setminus \{0\}$ is Lagrangian. 

The 1-form $\langle JH_{C(\Sigma)}, \cdot \rangle$ on the Lagrangian submanifold $C(\Sigma)\setminus\{0\}$ is closed. 
(This follows directly either from Theorem 3.4 of Dazord in \cite{HL}, or can be verified by locally exactness via the local expression $$
H_{C(\Sigma)}=-J\nabla_{C(\Sigma)}\theta
$$ 
given in \cite{HL}; this will be done in next section.)  
Its restriction along $\Sigma$ is therefore a closed 1-form $i^*\langle JH_{C(\Sigma)},\cdot\rangle$, as the pullback by the inclusion $i:\Sigma\to C(\Sigma)$ of a closed 1-form. Since the first Betti number of $\Sigma$ is zero 
($g_\Sigma = 0$), there is a smooth function $\theta_\Sigma$ on $\Sigma$ such that 
$$
d\theta_\Sigma = i^*\langle JH_{C(\Sigma)},\cdot\rangle.
$$  
Then
$$
\langle \nabla_\Sigma\theta_\Sigma, \cdot \rangle = d\theta_\Sigma=\langle JH_\Sigma, \cdot \rangle.
$$


As we have seen, the Hamiltonian stationary condition on $C(\Sigma)$ implies 
\begin{eqnarray*}
0&=& \operatorname{div}_{\Sigma} (JH_{\Sigma})\\
&=& \operatorname{div}_{\Sigma}(\nabla_\Sigma\theta_\Sigma)\\
&=&\Delta_g\theta_\Sigma.
\end{eqnarray*}

On the closed surface $\Sigma,$ then $\theta_\Sigma$ is constant, in turn, $\Sigma$
is minimal.

An immersed minimal Legendrian 2-sphere in ${\mathbb S}^5$ is totally geodesic. This is a known fact, for a proof, see for example in \cite{CY}.

(2) $g_{\Sigma}\geq1:$ As in the $g_{\Sigma}=0$ case, the isotropic and
Hamiltonian stationary condition gives us a local holomorphic function
$\left\langle Ju,u_{z}\right\rangle $ and global holomorphic Hopf 
1-differential $\left\langle Ju,u_{z}\right\rangle dz$. We only consider the
case that $\left\langle Ju,u_{z}\right\rangle dz$ is not identically $0.$ The
zeros of $\langle Ju,u_{z}\rangle$ are therefore isolated and near each of the
zeros, we can write
\[
\left\langle Ju,u_{z}\right\rangle =h\left(  z\right)  z^{k}%
\]
where $h$ is a local holomorphic function non-vanishing at the zero point
$z=0$ and $k$ is a positive integer. One can also view
\[
\left\langle Ju,u_{z}\right\rangle =\frac{1}{2}\left(  \left\langle
Ju,u_{1}\right\rangle -\sqrt{-1}\left\langle Ju,u_{2}\right\rangle \right)
\]
as the tangent vector
\[
\frac{1}{2}\left\langle Ju,u_{1}\right\rangle u_{1}-\frac{1}{2}\left\langle
Ju,u_{2}\right\rangle u_{2}=\frac{1}{2}\left\langle Ju,u_{1}\right\rangle
\partial_{1}-\frac{1}{2}\left\langle Ju,u_{2}\right\rangle \partial_{2}%
\]
along the tangent space $T\Sigma,$ where $\partial_{i}=\partial u/\partial
{t^{i}}$. The projection $\Pr Ju$ on the tangent space of $T\Sigma$ is locally
represented as
\[
\Pr Ju=\frac{\left\langle Ju,u_{1}\right\rangle \partial_{1}+\left\langle
Ju,u_{2}\right\rangle \partial_{2}}{\varphi^{2}}.
\]
The index of the globally defined vector field $\Pr Ju$ at each of its
singular points, i.e. where $\Pr Ju=0$, is the negative of that for the vector field
$\frac{1}{2}\left\langle Ju,u_{1}\right\rangle \partial_{1}-\frac{1}%
{2}\left\langle Ju,u_{2}\right\rangle \partial_{2}.$ Note that the index of
the latter is $k.$

From the Poincar\'{e}-Hopf index theorem, for any vector field $V$ with
isolated singularities on $\Sigma$, one has
\[
\sum_{V=0}\mbox{index}\left(  V\right)  =\chi\left(  \Sigma\right)
=2-2g_{\Sigma}\leq0.
\]

The zeros of $\Pr Ju$ are just the Legendrian points on $\Sigma$. So we
conclude that the number of Legendrian points is $2g_{\Sigma}-2$ counting the
multiplicity. This completes the proof of Theorem 1.1.

\begin{remark}
\label{remark1} \emph{As mentioned in the introduction, Theorem \ref{thm-Hopf}
(except the totally geodesic part) and its generalization to higher
codimensions can be extended for the non-orientable links. This can be seen as
follows. The Poincar\'{e}-Hopf index theorem holds on compact non-orientable
surfaces, our count of the indices of the still globally defined $\Pr Ju$ via
\emph{local} holomorphic functions is valid too, and the index of a singular
point of a vector field is independent of local orientations. Moreover, this
index counting argument yields an alternative proof for Theorem \ref{thm-Hopf}
(1) (except the totally geodesic part) and its generalization. }
\end{remark}

\section{Harmonic forms and Proof of Theorem \ref{thm-Hodge}}

Consider an immersed isotropic Hamiltonian stationary submanifold in
$S^{2n-1}$
\[
u:L^{m}\rightarrow S^{2n-1}\subset\mathbb{C}^{n}.
\]
The isotropy condition reads for any local coordinates $\left(  t^{1}%
,\cdots,t^{m}\right)  $ on $L^{m}$
\[
\langle Ju_{i},u_{j}\rangle =0
\]
with $J$ being the complex structure of $\mathbb{C}^{n}$ and $u_{i}=\partial
u/\partial t^{i}.$ 

The Hamiltonian stationary condition for the cone $C\left(
\Sigma\right)  =ru\left(  t\right)  $ is
\begin{eqnarray*}
0  &  =&\operatorname{div}_{C\left(  L\right)  }\left(  JH_{C\left(  L\right)
}\right) \\
& =& \left\langle \partial_{r}\left(  JH_{C\left(  L\right)  }\right)
,\partial_{r}\right\rangle +\operatorname{div}_{L}\left(  J\left(\frac{1}
{r}H_{L}\right)\right) \\
&  =&-\frac{1}{r^{2}}\left\langle JH_L,u\right\rangle +\frac{1}{r^{3}
}\operatorname{div}_{L}\left(  JH_{L}\right).
\end{eqnarray*}
Notice that $\langle JH_L,u\rangle$ and $\operatorname{div}_{L}\left(  JH_{L}\right)$ are independent of $r$. Therefore, the  equation above splits into two equations
\[
\operatorname{div}_{L}\left(  JH_{L}\right)  =0
\]
and%
\[
0=\langle JH_{L},u\rangle =-\langle \Delta_{g}u,Ju\rangle ,
\]
where $g$ is the induced metric on $L$ and $\Delta_g$ is the Laplace-Beltrami operator of $(L,g)$. 

To measure the deviation of the corresponding cone $C(u(L^{m}))$ from
being isotropic, we project $Ju$ onto the tangent space of $u(L^{m})$ in
$\mathbb{C}^{n}$. Note that the projection is the vector field along $u(L)$
\[
\hbox{Pr}Ju=\sum_{i,j=1}^{m}g^{ij}\langle Ju,u_{i}\rangle u_{j}%
\]
where $g_{ij}=\langle u_{i},u_{j}\rangle,1\leq i,j\leq m.$ The corresponding
1-form
\[
\alpha=\sum_{i=1}^{m}\left\langle Ju,u_{i}\right\rangle dt^{i}%
\]
is of course globally defined on $L^{m}$. In fact it is a harmonic 1-form,
because $\alpha$ is closed and co-closed as verified as follows: 
\begin{eqnarray*}
d\alpha &  =&\sum_{i,j=1}^{m}\left\langle Ju,u_{i}\right\rangle _{j}%
dt^{j}\wedge dt^{i}\\
&=&\sum_{i,j=1}^{m}\left(\langle Ju_{j},u_{i}\rangle +\langle Ju,u_{ij}\rangle \right)  dt^{j}\wedge
dt^{i}\\
&  =&\sum_{i,j=1}^{m}\langle Ju,u_{ij}\rangle dt^{j}\wedge dt^{i}\\
&=&0,
\end{eqnarray*}
and 
\begin{eqnarray*}
\delta\alpha &  =&\left(  -1\right)  ^{m\cdot1+m+1}\ast d\ast\alpha\\
&  =&-\ast d\left(  \sum_{i,j=1}^{m}\left(  -1\right)  ^{j+1}\sqrt{g}%
g^{ij}\left\langle Ju,u_{i}\right\rangle dt^{1}\wedge\cdots\wedge
\widehat{dt^{j}}\wedge\cdots\wedge dt^{m}\right) \\
&  =&-\ast\sum_{i,j=1}^{m}\partial_{j}\left(  \sqrt{g}g^{ij}\left\langle
Ju,u_{i}\right\rangle \right)  dt^{1}\wedge\cdots\wedge dt^{j}\wedge
\cdots\wedge dt^{m}\\
&  =&-\frac{1}{\sqrt{g}}\sum_{i,j=1}^{m}\partial_{j}\left(  \sqrt{g}%
g^{ij}\langle Ju,u_{i}\rangle \right) \\
&  =&-\sum_{i,j=1}^{m}\left( \left\langle Ju_{j},g^{ij}u_{i}\right\rangle
+\left\langle Ju,\frac{1}{\sqrt{g}}\partial_{j}\left(  \sqrt{g}g^{ij}%
u_{i}\right) \right\rangle \right) \\
&  =&-\left\langle Ju,\Delta_{g}u\right\rangle\\
& =&0,
\end{eqnarray*}
where we have used the isotropy condition and the consequence of Hamiltonian stationary condition in the last two steps, respectively.

The Hodge-de Rham theorem implies that the harmonic 1-form $\alpha$ must
vanish because the first Betti number of $L^m$ is zero by assumption. It follows that
$\hbox{Pr}Ju$ must vanish. Therefore, the cone $C(L^{m})$ is isotropic.

Next, we claim that the differential 1-form 
$$
\beta= \langle JH_L,\cdot\rangle
$$ 
on $L^m$ is closed. When $m=n-1,$ the isotropic cone $C(L^{n-1})$ is Lagrangian. By \cite{HL}, around each point of $C(L^{n-1})\setminus\{0\}$, there is a locally  defined Lagrangian angle $\theta$ such that
\[
H_{C\left(  L\right)  }= -J\nabla_{C\left(  L\right)  }\theta.
\]
Now the globally defined 1-form $\beta$ on the link $L$ can be expressed locally as 
$$
\beta = \langle \nabla_{C(L)}\theta,\cdot\rangle
=\langle\nabla_L\theta,\cdot\rangle = d_L\theta
$$
by noticing that $H_{C(L)} = H_L$ as $r=1$, where the second equality holds as the two 1-forms are on $TL$ and the tangent vectors to $L$ are orthogonal to $\partial_r$, and $d_L$ stands for the exterior differentiation on $L$. We conclude that $\beta$ is a closed 1-form on $L$. When $m <n-1$, the 1-form $\langle JH_{C(L)},\cdot\rangle$ is closed by assumption, so its restriction $\beta$ on $L$ is closed. 

Since the first Betti number of $L$ is zero, there is a smooth function 
$\theta_L$ on $L$ such that 
$$
\langle JH_L,\cdot \rangle = d_L\theta_L.
$$
This implies that the projection of $JH_L$ onto $TL$ satisfies 
$$
\sum_{i=1}^m \langle JH_L, E_i\rangle E_i=\nabla_L\theta_L
$$
where $\{E_1, \dots, E_m\}$ is a local orthonomal frame of $TL$. 
The Hamiltonian stationary condition on $C(L)$ asserts, as we have seen earlier, 
that 
\begin{eqnarray*}
\Delta_L\theta_L&=& \operatorname{div}_L\nabla_L\theta_L\\
&=&\sum^m_{i=1}\operatorname{div}_L \left(\langle JH_L,E_i\rangle E_i \right) \\
&=& \sum^m_{i=1} \left\langle\nabla_L \langle JH_L,E_i\rangle, E_i\right\rangle +\langle JH_L,E_i\rangle \operatorname{div}_L E_i\\
&=&\operatorname{div}_L (JH_L)+\langle JH_L,H_L\rangle\\
&=& 0.
\end{eqnarray*}
On the closed submanifold $L,$ then $\theta_L$ is constant.
In turn, for $m=n-1,$ $C\left(  L^{n-1}\right)  $ is minimal, and for $m<n-1,$
$C\left(  L^{m}\right)  $ is partially minimal, namely $H_{C\left(
L^{m}\right)  }$ vanishes on the normal subbundle $JTC\left(  L^{m}\right)  .$
The proof of Theorem 1.2 is complete.

\begin{remark}
\label{remark2} \emph{As the projection }$\Pr Ju$\emph{ and the adjoint
operator }$\delta$\emph{ is independent of the local orientations and the
Hodge-de Rham theorem holds for compact non-orientable manifolds (cf. [LW,
p.125--126]), we see that Theorem \ref{thm-Hodge} remains true for
non-orientable links }$L^{m}.$
\end{remark}

\begin{remark}\emph{
For a surface link $L^2\subset {\mathbb S}^{2n-1}$ with $g_L=0$ for the case $n>3$, if it is isotropic and $C(L^2)$ 
is Hamiltonian stationary, the same argument as in \cite{CY} leads to the conclusion that the second fundamental form of $L$ in ${\mathbb S}^{2n-1}$ vanishes 
in the normal subbundle $Ju\oplus JTL$.  When $n=3$, $L$ is totally geodesic in $\mathbb S^5$ as noted before.}
\end{remark}

\begin{corollary}
\label{cor-geo}Let $L^{m}$ be a compact immersed isotropic submanifold in the
unit sphere $S^{2n-1}\subset\mathbb{C}^{n}$. If the Ricci curvature of $L^{m}$
is nonnegative, and it is positive somewhere or the Euler characteristic
$\chi(L^{m})$ is not zero, then the Hamiltonian stationary cone $C(L^{m})$ is
isotropic, in particular, $C(L^{n-1})$ is Lagrangian (or equivalently
$L^{n-1}$ is Legendrian) and minimal when $m$ is the top dimension $n-1.$
\end{corollary}

Under the above condition, from Bochner [B, p.381], it follows immediately
that the first Betti number of $L^{m}$ is zero. Then Theorem 1.2 and its
non-orientable version imply the corollary.

\end{document}